\documentclass{article}
\usepackage{amssymb}
\usepackage{amsmath,amsfonts,amsthm} 
\newtheorem{alg}{Algorithm}
\newtheorem{definition}{Definition}
\newtheorem{theorem}{Theorem}
\newtheorem{example}{Example}
\newtheorem{remark}{Remark}
\newtheorem{proposition}{Proposition}

\usepackage{graphicx}
\def\T{{\rm T\,}}
\def\Td{{\rm Td\,}}

\def\Can{{\rm Can}}

\begin{document}

\title{Gr\"obner bases and combinatorics for binary codes.%\thanks{Grants or other notes
}


\author{M. Borges-Quintana\and M.A.~Borges-Trenard\thanks{Dpto. de Matem\'atica, FCMC, U. de Oriente, Santiago de Cuba, Cuba.
	      Partially supported by Spanish MEC grant SB 2003-0287.
              mijail@mbq.uo.edu.cu, mborges@mabt.uo.edu.cu }\and P.~Fitzpatrick\thanks{Boole Centre for Research in Informatics, UCC, Cork, Ireland. fitzpat@ucc.ie} \and E.~Mart\'{\i}nez-Moro\thanks{Dpto. de Matem\'atica Aplicada, U. de Valladolid, Castilla, Spain.
	      Partially supported by Spanish MEC MTM2004-00876 and MTM2004-00958.edgar@maf.uva.es}
}

%\date{Received: date / Revised: date}

\maketitle
\begin{abstract} 

In this paper we introduce a binomial ideal derived from a binary
linear code. We present some applications of a Gr\"obner basis of
this ideal with respect to a total degree ordering. In the first
application we give a decoding method for the code. By associating
the code with the set of cycles in a graph, we can solve the
problem of finding all codewords of minimal length (minimal cycles
in a graph), and show how to find a minimal cycle basis. Finally
we discuss some results  on the computation of the Gr\"obner
basis.

%\keywords{Binary codes \and Cycle bases \and  Gr\"obner bases.}
\end{abstract}

\section{Introduction}\label{int}

We associate with a binary linear code a Gr\"obner basis for total
degree compatible orderings such as degrevlex (Degree Reverse
Lexicographic), for which Gr\"obner bases are known to be easier
to compute. In our particular application the Gr\"obner basis has
additional properties that allow us to formulate an 
algorithm, which has the flavour of an FGLM approach and it is especially adapted to our setting. We show how
the Gr\"obner basis of the code can be used for decoding and solve
several problems related to graphs associated with the code.

In the paper we use the term {\em code} to refer only to {\em
binary linear code} even though some of our results ({\em cf.\/}
Sections~\ref{s:bc-Mon},\ref{s:id-bc}) can be extended to the
non-binary case (see also \cite{BBW,bbw-rep,bbm}). The outline of
the paper is as follows. In Section~\ref{s:bc-Mon} we define a
monoid connected to a binary linear code. An ideal associated with
the code is introduced in Section~\ref{s:id-bc} together with a
decoding method that makes use of a Gr\"obner basis of the ideal.
In fact, decoding is carried out using classical reduction to the
canonical form. Some other applications are developed in
Section~\ref{s:app} such as finding all minimal cycles in a graph
and a minimal cycle basis. In Section~\ref{s:comp} a linear
algebra procedure (related to FGLM) is used to compute the
Gr\"obner basis for the ideal associated with a code. This method
is applicable in a general setting but in our setting it has
additional computational advantages.

\section{Binary linear codes and monoids}\label{s:bc-Mon}
\subsection{Binary linear codes}\label{sb:lc}
  Let
$\mathbb{F}_2$ be  the finite field with $2$ elements. A linear
code $\mathcal C$ of dimension $k$ and length $n$  is the image of
a linear mapping $L: \mathbb{F}_2^k\to \mathbb{F}_2^n$,  where
$k\leq n$, i.e. $\mathcal{C}=L(\mathbb{F}_2^k)$. There exists a
$n\times (n-k)$  matrix $H$, called a {\sl parity check matrix},
such that $cH=0\,$ if and only if $c\in\mathcal C$. On the other
hand, there exists a $k\times n$ generator matrix $G$ such that
$\mathcal{C}=\{u\,G\;\mid\;u\in \mathbb{F}_2^k\}$. Normally, we
consider check matrices to have linearly independent columns and
generator matrices linearly independent rows. However, in some
situations it is useful to regard as a check matrix any matrix
whose left nullspace is the code, and as generator matrix any
matrix whose row space is the code. The weight of a codeword is
its Hamming distance to the word $0$, and the minimum distance $d$
of a code is the minimum weight among all the non-zero codewords.
The error correcting capacity of a code is
$t=\left[\frac{d-1}{2}\right]$, where $\left[\cdot\right]$ is the
greatest integer function. Let $B({\mathcal C},t)=\{y\in
\mathbb{F}_2^n\mid \exists c\in\mathcal{C}\hbox{ s.t. } d(c,y)\leq
t \} $, it is well known that the equation \[eH=yH\] has a unique
solution $e$ with $\mathrm{weight}(e)\leq t$ for $y\in B({\mathcal
C},t)$.
\subsection{The monoid associated with a binary code}\label{sb:mon-bc}
Let $[X]$ be the free commutative monoid  generated by the $n$
variables $X=\{ x_{1},\dots ,  x_{n}\}$. We have the following map
from $X$ to $ \mathbb{F}_2^n$:
\begin{equation}\begin{split}
\psi: X\to & \mathbb{F}_2^n\\
x_{i}\mapsto \,& \,{\mathbf e}_i=(0,\dots,0,\underbrace{1}_{i},0,\dots ,0)
\end{split}\end{equation}
The map $\psi$ can be extended to a morphism from  $[X]$ onto
$\mathbb{F}_2^n$, where
\begin{equation}\psi\left( \prod_{i=1}^n
x_{i}^{\beta_{i}}\right)=\left(\beta_{1}\bmod 2 ,\dots,
\beta_{n}\bmod 2 \right)\end{equation} When no confusion arise we
will use $x_i$ to refer the indeterminate in the monoid or the
associated vector $\mathbf{e}_i$ in $\mathbb{F}_2^n$. A code
$\mathcal C$ defines an equivalence relation $R_{\mathcal C}$ in
$\mathbb{F}_2^n$ given by
\begin{equation}\label{eq:eq-rel}
(x,y)\in R_{\mathcal C}\Leftrightarrow x-y\in \mathcal C.
\end{equation}
If we define $\xi (u)=\psi (u)H$, where $u\in [X]$, the above
congruence can be translated to $[X]$ by the morphism $\psi$ as
\begin{equation}u\equiv_{\mathcal C} w\Leftrightarrow (\psi(u),\psi(w))\in R_{\mathcal
C}\Leftrightarrow \xi(u)=\xi(w). \end{equation}
 The morphism $\xi$ represents the
transition of the syndromes from $\mathbb{F}_2^n$ to $[X]$. Thus,
 $\xi(w)$ is the syndrome of $w$, which is the
syndrome of $\psi(w)$.

For the sake of simplicity we will use $w,u,v$ as words in $[X]$
and vectors in $\mathbb{F}_2^n$, as long as the meaning is clear
from the context. The connection between the two structures can be
understood from the following setting
\begin{equation} w=1\cdot x_{i_1}\cdot \ldots \cdot x_{i_m}
\in [X]\rightarrow \psi(w)=\mathbf{0}+\mathbf{e}_{i_1}+\ldots+\mathbf{e}_{i_m}
\in \mathbb{F}_2^n.\end{equation}

\begin{definition}[standard word]\label{st-rep} The word $w=\prod_{i=1}^n\,
x_{i}^{\beta_{i}}$ is said to be {\em standard\/} if $\beta_{i} <
2$, for every $i \in \{1,\ldots,n\}$. Given $y\in {\mathbb F}_2^n$
we say that $w$ is the {\em standard representation\/} of $y$ if
$\psi(w)=y$ and $w$ is standard .
\end{definition}

\subsection{Binary codes and the set of cycles in a graph}\label{sb:bib-graf}
Let $G=(V,E)$ an undirected 2-connected graph without loops or
multiple edges, where $V$ is the set of vertices and $E$ the set
of edges. An edge is denoted by an unordered pair of vertices
$(x,y)$. A cycle is a subgraph such that any vertex degree is
even. Therefore, a cycle can be written as either a set of edges
$\{(x_1,x_2),(x_2,x_3),\ldots\}$ or as a closed path
$(x_1,x_2,x_3,\ldots,x_1)$. The length of a cycle is the number of
edges it contains.

The sum of two cycles is defined as the symmetric set difference
$C+C^\prime=(C\cup C^\prime)\setminus (C\cap C^\prime)$. With this
sum the set of cycles forms an $\mathbb{F}_2$-vector space which
is a subspace of $\mathbb{F}_2^m$, where $m=| E |$ (the number of
edges). Therefore, the set $\mathcal{C}$ of cycles in a graph can
be considered as a binary code of length $m$. A basis of this
vector space is called a {\em cycle basis}, and its dimension is
well-known to be the Betti number $dim(\mathcal{C})=m - |V| +1$
(see, for example, \cite{Leydolda,pewe,vismara}). We define the
{\em length of a basis\/} as the total length of the cycles in it.

\section{The ideal associated with a code}\label{s:id-bc}
In this section we define a particular ideal associated with a
code. Ideals associated with soft-decision maximum likelihood
decoding can be found in \cite{ikegami}, and these in turn are
related to ideals arising in integer programming using Gr\"obner
basis \cite{clo}.

Consider the polynomial ring $K[X]$, where $K$ is a field. Let
$<=<_\T$ be a fixed, total degree compatible term order with $x_1<
x_2 \cdots < x_n$ on $[X]$. We use $<$ for this term order as the
meaning of the symbol will always be clear from the context, and
write $>$ where appropriate. As usual, $\T(f)$ denotes the
maximal term of a polynomial $f$ with respect to the order $<$ and
$\Td(f)$ the total degree of the maximal term $\T(f)$ of $f$. The
set of maximal terms of the set $F\subseteq K[X]$ is denoted
$\T\{F\}$ and $\T(F)$ denotes the semigroup ideal generated by
$\T\{F\}$. Finally, $\langle F \rangle$ is the polynomial ideal
generated by $F$.

\begin{definition}\label{d:id-cod} Let $\mathcal C$
be a code and  $R_{\mathcal C}$ the equivalence relation defined 
in equation (\ref{eq:eq-rel}). The ideal $I({\mathcal C})$
associated with $\mathcal C$ is
\begin{equation} I(\mathcal C)=\langle \{w -v \,\mid\, (\psi(w),\psi(u))\in R_{\mathcal C}
\} \rangle\subseteq K[X].
\end{equation}
\end{definition}

Let be $\{w_1,\ldots,w_k\}$ be the row vectors of a generator
matrix for a code (more generally any matrix whose rows span the
code $\mathcal{C}$), i.e., a basis (spanning set) of the code as
subspace of $\mathbb{F}_2^n$. Let
\begin{equation}\label{eq:Iinit}
I=\langle \{w_1-1,\ldots,w_k-1\}\cup\{x_i^2-1\mid i=1,\ldots, n\}
\rangle\end{equation} be the ideal generated by the set of
binomials $\{w_1-1,\ldots,w_k-1\}\cup\{x_i^2-1 \mid i=1,\ldots,
n\}\in K[X]$. Since $\{w_1,\ldots,w_k\}$ generate $\mathcal C$ it
is clear that $I=I({\mathcal C})$.

\subsection{Error-correcting reduced Gr\"obner basis}\label{sb:err-GB}

Let $G_\T$ be the reduced Gr\"obner basis of the ideal
$I(\mathcal{C})$ with respect to $<$. Note that $G_\T$ can be
computed by Buchberger's algorithm starting with the initial set
$\{w_1-1,\ldots,w_k-1\}\cup\{x_i^2-1\mid i=1,\ldots, n\}$.
However, there are some computational advantages in this case. The
coefficient field is $\mathbb{F}_2$ (and therefore there is no
coefficient growth), and the maximal length of a word appearing in
the computation is $n$ (the binomials $x_i^2-1$ prevent the length
being greater than $n$). Thus the two principal disadvantages of
Gr\"ober basis computations are not valid for this case. In
addition, total degree compatible term orders are among the most
efficient for the computation of Gr\"obner bases.

Although the usual reduction could be carried out with the same
result, we introduce a special reduction in order to have a more
efficient process.

\begin{definition}[One step reduction]\label{redG} Reduction in one step ($\longrightarrow $)
using $G_T$ is defined
as follows. For any $w \in [X]$:
\begin{enumerate}
\item reduce $w$ to its standard form $w^\prime$ using the
relations $x_i^2 \longrightarrow  1$, for all $x_i \in X$. \item
reduce $w^\prime$ with respect to $G_T$ by the usual one step
reduction.
\end{enumerate}
\end{definition}

This reduction process is well defined since it is confluent and
noetherian. Thus, it will end after a finite number of one step
reductions with a unique irreducible element corresponding to the
starting element. Moreover, if we denote by $\Can(w,G_T)$ the
canonical form of $w$ with respect to $G_\T$ we have the following
result.

\begin{theorem}[Canonical forms of the vectors in $B(\mathcal{C},t)$]\label{t:DecGB}\label{dec:t-RB}
Let $\mathcal C$ be a code and let $G_\T$ be the reduced Gr\"ober
basis with respect to $<$. If $w \in [X]$ satisfies the condition
$\mathrm{weight}(\psi(\Can(w,G_\T))) \leq t$ then
$\psi(\Can(w,G_\T))$ is the error vector corresponding to
$\psi(w)$. On the other hand, if
$\mathrm{weight}(\psi(\Can(w,G_\T)))
> t$ then $\psi(w)$ contains more than $t$ errors.
\end{theorem}
\begin{proof}
The uniqueness of the canonical form is guaranteed by its
definition, and thus we only need to prove that the standard
representation of the error vector associated with a vector $y$
satisfies the condition to be the canonical form of $y$.

Let $y \in B(\mathcal{C},t)$ and denote by $e=e_y$ be the error vector
corresponding to $y$. Then $eH = yH$ and $\mathrm{weight}(e) \leq
t$. If $w_{e}$ is the standard representation of $e$ then
$\mathrm{weight}(e)$ coincides with the total degree of $w_{e}$.
Accordingly, $\Td (w_\mathbf{e}) \leq t$. It is clear that there
cannot be another word $u$ such that $\Td (u)\leq t$ and
$\xi(u)=yH$, since this would mean that there are two solutions
for the linear system with weight at most $t$, and this is not
possible because $y \in B(\mathcal{C},t)$. Therefore, it is clear that $w_e$
is the minimal element with respect to $<$ having the same
syndrome as $y$.\qed
\end{proof}

We see later that the error-correcting capability $t$ of the code
can be computed from $G_\T$ (see Remark~\ref{rem:t}).

\begin{example}[Decoding a binary code using its associated Gr\"obner basis]\label{ej:dec_GB_bin}
Let be $G$ be a generator matrix of the $[6,2,3]$ binary code
$\mathcal{C}$ over $\mathbb{F}_2^6$ defined as

$$G=\left( \begin{array}{llllll}
1 & 1 & 0 & 1 & 1 & 0\\
0 & 1 & 1 & 0 & 0 & 1\\
1 & 0 & 1 & 1 & 1 & 1
\end{array}
 \right)$$
 From this matrix we obtain a set of generating polynomials for $I(\mathcal{C})$
 as in equation (\ref{eq:Iinit}) as follows:
 \begin{equation*}\begin{split}
I(\mathcal{C})= & \langle x_1x_2x_4x_5-1,x_2x_3x_6-1,x_1x_3x_4x_5x_6-1, \\
& x_1^2-1,x_2^2-1,x_3^2-1,x_4^2-1,x_5^2-1,x_6^2-1 \rangle.\end{split}
\end{equation*}

A Gr\"obner basis of $I(\mathcal{C})$ with respect to the
degrevlex is
\begin{equation*}\begin{split}
G_\T= & \{x_1^2-1,x_2^2-1,x_3^2-1,x_4^2-1,x_5^2-1,x_6^2-1,\\  &
x_2x_3 -x_6, x_2x_4-x_1x_5, x_2x_5-x_1x_4, x_2x_6-x_3, \\ &
x_3x_6-x_2, x_4x_5-x_1x_2, x_1x_3x_4 -x_5x_6, x_1x_3x_5-x_6x_4, \\
& x_6x_1x_4-x_3x_5, x_1x_5x_6-x_3x_4\}. \end{split}
\end{equation*}
%where for convenience the elements $x_i^2-1$ have been rearange at the begining of the list.
The decoding process consists of obtaining the errors as a common
reduction process modulo the Gr\"obner basis $G_\T$. Suppose the
word $w=x_1x_2x_3x_4x_5$ is received. The canonical form of $w$
modulo $G_\T$ is $x_3$, since $x_3$ has weight 1 and the code is
1-error correcting (see Theorem~\ref{t:DecGB}), then the
corresponding codeword is $x_1x_2x_4x_5$ or $110110$.
\end{example}

\section{Further Applications.}\label{s:app}

We will show that the Gr\"ober basis $G_\T$ for a code $\mathcal{C}$  can be used to
solve some other problems in coding theory and graph theory.
In general, let $c_g$ be the codeword associated to the binomial
$g=w-v \in I(\mathcal{C})$, so that $c_g=\psi(w)+\psi(v)$, and let
$w_c$ be the standard word corresponding to $c$.

\begin{theorem}[Reduction of a codeword]\label{t:red-we}
Let $c$ be a codeword such
that $\mathrm{weight}(c)=d^\prime$. Then there exists $g_1 \in G_\T$ 
such that:
\begin{enumerate}
\item $\Td(g_1)\leq t^\prime=[(d^\prime -1)/2]+1$. 
\item $w_c\stackrel{g_1}{\longrightarrow} w_2$, such that $\Td(w_2)\leq
d^\prime $ and $w_2 < w_c$. 
\item $c=c_{g_1}+c_{w_2}$, where
$\mathrm{weight}(c_{g_1})\leq d^\prime $ and
$\mathrm{weight}(c_{w_2})\leq d^\prime $.
\end{enumerate}
\end{theorem}
\begin{proof}
Let $w_{c_1}$ and $u_{c_1}$ be such that $w_c=w_{c_1}u_{c_1}$ with
$\mathrm{weight}(w_{c_1})=t^\prime$ and $u_{c_1}<w_{c_1}$. It is
clear that $w_{c_1}-u_{c_1} \in I(\mathcal{C})$ since $w_{c_1}$
and $u_{c_1}$ have the same syndrome. Therefore $w_{c_1} \in
\T(G_\T)$. Let $g_1=w_1-v_1 \in G_\T$ where $v_1 < w_1$ satisfy $w_{c_1}=w_1u_1$ for
some $u_1 \in [X]$. Then $\Td(g_1)=\Td(w_1) \leq t^\prime$ and
hence $g_1$ satisfies condition (\textit{1.}).

Now, $w_c \stackrel{g_1}{\longrightarrow} w_2=u_{c_1}u_1v_1$. Note
that $\Td(u_{c_1})=d^\prime -t^\prime$ and $\Td(u_1v_1)\leq
t^\prime$, which implies $\Td(w_2)\leq d^\prime $. Thus,
$\mathrm{weight}(w_2)\leq d^\prime $. Also $v_1 < w_1$ implies
$v_1u_1u_{c_1}< w_1u_1u_{c_1}$, that is, $w_2< w_c$, and (\textit{2}.)
follows.

In order to prove (\textit{3}.), we observe that $\Td(w_1)\leq \Td(w_{c_1})$
and by the construction of $w_{c_1}$ and $u_{c_1}$ and $v_1$ being a canonical form, we have also
that $\Td(v_1)\leq \Td(u_{c_1})$. Thus, $\Td(w_1v_1)\leq d^\prime
$ and $\mathrm{weight}(c_{g_1})\leq d^\prime $. It is easy to see
that $c=c_{g_1}+c_{w_2}$ since
$c=\psi(w_1)+\psi(u_1)+\psi(u_{c_1})$,
$c_{g_1}=\psi(w_1)+\psi(v_1)$, and
$c_{w_2}=\psi(u_{c_1})+\psi(u_1)+\psi(v_1)$.\qed
\end{proof}
\begin{remark}\label{rem:t} As a consequence
 of (\textit{1}.) in this Theorem, given $G_\T$ and $g\in G_\T$ a
binomial such that $$\Td(g)=\min\left\lbrace \Td(f) \mid f \in
G_\T\setminus \{x_i^2-1 \mid i=1,...,n\} \right\rbrace $$ we have
$t=\Td(g)-1$. Moreover, in order to find such a $g$ it is not
necessary to compute the whole Gr\"obner basis $G_\T$ (see
Theorem~5 in \cite{bbw-rep} or Remark~\ref{rem:tpat} in
Section~\ref{s:comp}).
\end{remark}

The following propositions provides important properties of
$\T\{G_\T\}$.
\begin{proposition}(Relation between $\Td(g)$ and $\mathrm{weight}(c_g)$)\label{p:heads}
Let $g\in G_\T$ satisfy $\mathrm{weight}(c_g )=d^\prime$, and let
$t^\prime=[(d^\prime-1)/2]+1$. Then $\Td(g)=t^\prime$ or
$\Td(g)=t^\prime+1$.
\end{proposition}
\begin{proof}
It is clear that if $t^\prime=1$ (which would imply that the code
$\mathcal{C}$ has 0 error-correcting capability) the result is
true. We may assume that $t^\prime\geq 2$.

Obviously $\Td(g)\geq t^\prime$, otherwise
$\mathrm{weight}(c_g)<d^\prime$. Suppose that $\Td(g) >
t^\prime+1$, and let $\T(g)=xw, \Can(g,G_T)=v$ (where $x$ is any
variable belonging to the support of $\T(g)$). Observe that
$\Td(w)\geq t^\prime+1$ and $\Td(v)\leq t^\prime-2$. As a
consequence, $w > xv$ and thus $w \in \T(G_\T\setminus \{g\})$
(note that $w-xv\in I(\mathcal{C})$) which cannot happen because
$G_\T$ is a reduced Gr\"obner basis. This completes the proof.\qed
\end{proof}

\begin{proposition}(Codewords of minimal weight)\label{p:minpes}
Let $c$ be a codeword of minimal weight $d$. If $d$ is odd then
there exists $g \in G_\T$ such that $c=c_g$ and $\Td(g)=t+1$. If
$d$ is even then either there exists $g \in G_\T$ such that
$c=c_g$ and $\Td(g)=t+1$ or there exist $g_1,\, g_2\in G_\T$ such
that $c=c_{g_1}+c_{g_2}=\psi(w_1)+\psi(w_2)$, where $g_1=w_1-v$,
$g_2=w_2-v$
$(w_1=\T(g_1),w_2=\T(g_2),\,v=\Can(g_1,G_\T)=\Can(g_2,G_\T))$,
with $t+1=\Td(g_1)=\Td(g_2)$.
\end{proposition}
\begin{proof}
Let $f=w_{c}-v_{c}$, where $\Td(w_{c})=t+1$, $\Td(v_{c})=d-t-1$ and
$c_f=c$.

If $d$ is odd then, by Theorem~\ref{t:DecGB}, $v_c$ is a canonical
form ($\mathrm{weight}(v_c)=t$) and $\Can(w_c,G_\T)=v_c$. By
Theorem~\ref{t:red-we}, there exists $g_1\in G$ that satisfies the
conditions of the theorem. In this case, part (\textit{1}.) implies that
$\Td(g_1)=t+1$. (By Proposition~\ref{p:heads} there are no maximal
terms of degree less than $t+1$, apart from the monomials with
support size $1$). Consequently, $T(g_1)=w_c$ and therefore,
$f=g_1$.

If $d$ is even then $\mathrm{weight}(v_c)=t+1$ and it is not
necessarily a canonical form. If it is a canonical form then we
are in the same case as before, that is, there exists $g_1 \in G$
such that $c=c_{g_1}$ and $\Td(g_1)=t+1$. If $v_c$ is not a
canonical form then there exist $g_1,g_2 \in G$, such that
$\Td(g_1)=\Td(g_2)=t+1$, $\T(g_1)=w_c$, $\T(g_2)=v_c$, and
$\Can(g_1,G_\T)=\Can(g_2,G_\T)=v$. It is easy to check that these
two binomials satisfy $c=c_{g_1}+c_{g_2}$
($c=\psi(w_c)+\psi(v_c)=c_{g_1}+c_{g_2}$ because the term
$\psi(v)$ appears twice and therefore vanishes).\qed
\end{proof}

Using the connection between cycles in graph and binary codes (see
Section~\ref{sb:bib-graf}), the previous theorem enables us to
obtain all the minimal cycles of a graph according to their
lengths. We will use $G_\T$ to compute a minimal cycle basis (see
\cite{Leydolda}), that is, a basis of the set of cycles considered
as vector space which has minimal length. First, we have the
following result, whose proof is a straightforward application of
Theorem~\ref{t:red-we}.

\begin{proposition}(Decomposition of a codeword)\label{p:decomp_cod}
Any codeword (or cycle in the corresponding graph) can be
decomposed as a sum of the form $c=\sum_{i=1}^l\,c_{g_i}$, where
$g_i \in G_\T$, $\mathrm{weight}(c_{g_i})\leq \mathrm{weight}(c)$,
and $$\Td(g_i)\leq \left[ \frac{(\mathrm{weight}(c)-1)}{2}\right] +1, \hbox{ for all
} i=1,\ldots, l.$$
\end{proposition}

By Theorem~\ref{t:red-we}, $c \in C$ can be reduced in one step
while the weight of $c_{g_1}$ and $c_{w_2}$ remains less than or
equal to $\mathrm{weight}(c)$. It is sufficient to carry this out
finitely many times because the reduction process must arrive at
the canonical form 1 (the empty word) after finitely many steps
($c_{emptyword}=(0,\ldots,0)$).

A minimal cycle basis can be obtained as a certain subset
$G^\prime$ of $G_\T$. The computation of $G_\T$ guarantees steps
similar to those in Horton's Algorithm for computing a minimal
cycle basis (see \cite{vismara}). A greedy algorithm can be used
to extract a cycle basis from the set $\{c_g\,\mid \, g \in
G_\T\}\setminus (0,\ldots,0)$, which turns out to be a minimal
cycle basis. This is made explicit in the following theorem.

\begin{theorem}[Finding a minimal cycle basis]\label{t:minbas}
Given the set $\mathcal{C}^\prime=\{c_g\,\mid \, g \in
G_\T\}\setminus (0,\ldots,0)$, where the elements of
$\mathcal{C}^\prime$ are ordered so that $c_{g_1} \prec c_{g_2}$
when one of the following conditions holds:
\begin{enumerate}
\item $\Td(g_1)<\Td(g_2)$. \item $\Td(g_1)=\Td(g_2)$ and
$\mathrm{weight}(c_{g_1})<\mathrm{weight}(c_{g_2})$. \item
$\Td(g_1)=\Td(g_2)$,
$\mathrm{weight}(c_{g_1})=\mathrm{weight}(c_{g_2})$, and $g_1 <
g_2$.
\end{enumerate}
Then the cycle basis obtained by applying a greedy algorithm
to $\mathcal{C}^\prime$ is a minimal cycle basis.
\end{theorem}

\begin{remark}
When $G_\T$ is computed it is close to being ordered according to
$\prec$. The only changes necessary are to reorder elements of the
same maximal term degree, by considering first the weights of the
corresponding codewords.
\end{remark}
\begin{proof}
There are two things to show in order to prove the result.
\begin{enumerate}
\item The set $\mathcal{C}^\prime$ contains a minimal cycle
basis. \item The ordering $\prec$ used to order the set
$\mathcal{C}^\prime$ is weight compatible with the goal of
obtaining a basis of minimal length.
\end{enumerate}
If these conditions hold then it is clear that a minimal cycle
basis will be obtained by applying a greedy algorithm to extract a
basis from $\mathcal{C}^\prime$. Since the set
$\mathcal{C}^\prime$ is a generating set of $\mathcal{C}$, it does
contain a basis.

Proof of ({1}.): Let $B=\{c_1,\ldots,c_l\}$ be a minimal cycle
basis. By applying Proposition~\ref{p:decomp_cod} we can decompose
any $c_i$ as
$$c_i=\sum_{j=1}^{n_i}\,c_{g_{ij}}, \hbox{ where }
\mathrm{weight}(c_{g_{ij}})\leq \mathrm{weight}(c_i) \hbox{ for
all }  j=1,\ldots,n_i.$$ Let $\mathcal{C}(B)=\{c_{g_{ij}}\mid
i=1,\ldots,l; j=1,\ldots,n_i\}$. Is clear that $\mathcal{C}(B)$ is
a generating set of $\mathcal{C}$. Moreover, the basis $B^\prime$
obtained by applying a greedy algorithm to $\mathcal{C}(B)$ has
length at most the length of $B$. Thus, $B^\prime$ is a minimal
cycle basis. Note that $\mathcal{C}(B) \subseteq
\mathcal{C}^\prime$.

Proof of ({2}.): Let $g_1,g_2 \in G_\T$ satisfy
$d_1=\mathrm{weight}(c_{g_1})<\mathrm{weight}(c_{g_2})=d_2$. Let
$t_1=[(d_1^\prime-1)/2]+1$ and $t_2=[(d^\prime-1)/2]+1$, so that
$t_1 \leq t_2$. The only conflict between $\prec$ and the weights
occurs when $\Td(g_1)>\Td(g_2)$ and this is possible only if
$\Td(g_1)>t_1$ (due to Proposition~\ref{p:heads} and the
inequality $t_1 \leq t_2$). By Proposition~\ref{p:decomp_cod} we
can find a set $\{c_{f_i} \mid i=1,\ldots,l\}$ such that
$c_{g_1}=\sum_{i=1}^l\,c_{f_i}$, where $f_i \in G_\T$,
$\mathrm{weight}(c_{f_i})\leq d_1$ and $\Td(f_i)\leq t_1$, for all
$i=1,\ldots , l$. This means that, in this case, $c_{g_1}$ is
already a linear combination of elements in $\mathcal{C}^\prime$
that occur earlier according to $\prec$. When $\Td(g_1)\leq
\Td(g_2)$ (and $d_1 < d_2$) we have $c_{g_1}\prec c_{g_2}$. This
completes the proof.\qed
\end{proof}
\begin{example}\label{sb:exmp} Given a graph
 $(V,U)$ of five vertices $V=\{1,2,3,4,5\}$ and six edges
 $U=\{(1,2),(1,4),(1,5),(2,3),(3,4),(4,5)\}$, the corresponding vector
 space is $\mathbb{F}_2^6$ (the length of codewords is the number of edges).
 It is easy to form a check matrix $H$ (whose columns are not, in general,
 linearly independent). Then $c \in \mathbb{F}_2^6$ is a cycle if and only if
 $cH=0$. Each row of $H$ corresponds to the representation of one of the
 edges such that there are exactly two ones in the positions corresponding
 to the vertices of the edge, so the matrix is as follows\vspace{0.2cm}
$$H=\left(
\begin{array}{lllll}
             1 & 1 &0 &0 &0\\
             1 &0 &0 &1 &0\\
             1 &0 &0 &0 &1\\
             0 &1 &1 &0 &0\\
             0 &0 &1 &1 &0\\
             0 &0 &0 &1 &1
\end{array}
\right).$$ From this matrix one can compute a generator matrix
$G$, although for this example it is easy to see that there are
just three cycles, which are those of Example~\ref{ej:dec_GB_bin}.
Thus, the matrix $G$ of that example is a generator matrix and we
have already computed the Gr\"ober basis $G_\T$ for this code.

\begin{description}
\item[\bf Application of Theorem~\ref{t:red-we}] Let us consider the
codeword (i.e. the cycle) $w_c=x_1x_3x_4x_5x_6$, $d^\prime=5$, and
$t^\prime=3$. Then it is clear that $x_1x_3x_4-x_5x_6\in
I(\mathcal{C})$, so $x_1x_3x_4 \in T(G_\T)$. Observe that
$g=x_1x_3x_4 -x_5x_6\in G_\T$ and $c=c_g$, which means that $c$ is
reduced to $(0,\ldots,0)$ in one step by $c_g$. Note that
$\Td(g)=3$.

\item[\bf Application of Proposition~\ref{p:minpes}] In this case the
minimum distance is $d=3$, Then all codewords (cycles) of minimal
weight (minimal length) can be obtained as certain $c_g$ where
$g\in G_\T$. In this case there is just one, namely,
$w_c=x_2x_3x_6$ ($c=(0,1,1,0,0,1)$).

\item[\bf Application of Proposition~\ref{p:decomp_cod}] Let
$g_1=x_3x_6-x_2,\,g_2=x_2x_4-x_1x_5$, we have that $g_1,g_2 \in G_\T$,
$c=c_{g_1}+c_{g_2}$, and all the conditions for
$\mathrm{weight}(\cdot)$ and $\Td(\cdot)$ are satisfied.

\item[\bf Finding a minimal cycle basis] We observe that
\begin{equation*}\begin{split}
\mathcal{C}^\prime= &\{c_{x_2x_3 -x_6},c_{x_2x_6-x_3},c_{x_3x_6-x_2},c_{x_2x_4-x_1x_5}, c_{x_2x_5-x_1x_4},c_{x_4x_5-x_1x_2}, \\ 
& c_{x_1x_3x_4-x_5x_6},c_{x_1x_3x_5-x_6x_4},c_{x_6x_1x_4-x_3x_5},c_{x_1x_5x_6-x_3x_4}\}
\end{split}\end{equation*}
where $\prec$ has been used to reorder the binomials at the same
level according to $\Td()$. Applying a greedy algorithm to
$\mathcal{C}^\prime$ we first choose $c_1=(0,1,1,0,0,1)$, the next
two binomials correspond also to $c_1$, and then the second
linearly independent vector, corresponding to $g=x_2x_4-x_1x_5$,
is $c_2=(1,1,0,1,1,0)$. Since the dimension of the vector space is
$2$, we already have a basis which is a minimal cycle basis by
Theorem~\ref{t:minbas}.
\end{description}
\end{example}

\section{Computation of the Gr\"obner basis.}\label{s:comp}

In this section we present a linear algebraic procedure that
allows us to compute the Gr\"obner basis associated with a code.
The background to this technique can be found in \cite{patrick1,patrick2}.

Given a set $F=\{ f_1, f_2, \ldots , f_r\}$ of polynomials in
$K[X]=K[x_1,\ldots ,x_n]$ generating an ideal $I$ let compute a
basis for the syzygy module $M$ in  $K[X]^{r+1}$ of the generator
set $F^\prime=\{-1,f_1, f_2, \ldots , f_r\}$. Each of the syzygies
corresponds to a solution
$$f=\sum_{i= 1}^r b_if_i \qquad b_i\in K[X], i=1,\ldots , r $$
and thus points to an element $f$ in the ideal $I$ generated by
$F$.

The main idea is that the set
\begin{equation}\begin{array}{c}
\mathbf{f}_1=(f_1, 1,0,0,\ldots , 0)\\
\mathbf{f}_2=(f_2, 0,1,0,\ldots , 0)\\
\vdots \\
\mathbf{f}_r=(f_r, 0,0,0,\ldots , 1)
\end{array}\end{equation}
is a basis of the syzygy module $M$, and moreover it is a
Gr\"obner basis with respect to a position over term (POT)ordering
$<_w$ induced from an ordering $<$ in $K[X]$ and the weight vector
$\mathbf{w}=(1,\mathrm{T}_<(f_1),\ldots ,\mathrm{T}_<(f_r) )$.
Also, the leading term  of $\mathbf{f}_i$ is $\mathbf{e}_{i+1}$
with respect the ordering $<_w$ where $\mathbf{e}_{j}$ denotes the
unit vector of length $r+1$ (see \cite{Adams} for an introduction
to Gr\"obner bases of modules).

Now we use the FGLM idea \cite{FGLM} and run through the terms of
$K[X]^{r+1}$ in the order determined by $<$ and $\mathbf{e}_i <
\mathbf{e}_j$ if $i<j$, using a term over position (TOP) ordering.
At each step the canonical form of the term with respect to the
original basis is $0$ apart from the first component so the
determination of the linear relations takes place in that
component. This provides a convenient representation for the
canonical form with respect to the initial Gr\"obner basis as a
$K$-vector space, and any linear relation obtained as a
consequence of reduction of the first component in $K[X]$ will
give a corresponding relation for the elements of the module.

\begin{example}
Let $I=\left\langle x^2+x+1,xy+x+1\right\rangle $ in $\mathbb
F_2[x,y]$ and take $<$ to be the deglex order with $x<y$.
Displaying only the first component we have

$$\begin{array}{l|cccccccccc}
 & 1 & x & y & x^2 & xy & y^2 & x^3 & x^2y & xy^2 & y^3\\ \hline
  (1,0,0) & 1 & &  & & & & & & & \\
  (0,1,0) & 1 & 1 &  & 1& & & & & & \\
  (0,0,1) & 1 & 1&  & &1 & & & & & \\ \hline
  \hbox{after reduction} &  & &  & & & & & & & \\
  (1,0,0) & 1 & &  & & & & & & & \\
  (1,1,0) &  & 1 &  & 1& & & & & & \\
  (0,1,1) &  & &  & 1& 1 & & & & & \\ \hline
    \hbox{introduce } x &  & &  & & & & & & & \\
  (x,0,0) &  & 1 &  &  & & & & & & \\
  (x,x,0) &  &  &  & 1& & & 1& & & \\
  (0,x,x) &  &  &  & &  & & 1 & 1&  & \\ \hline 
  \hbox{after reduction} &  & &  & & & & & & & \\
  (x+1,0,1) &  & &  & & 1& & & & & \\
  (1,x+1,0)  &  &  &  & & & & 1& & & \\
  (1,1,x) &  &  &  & &  & & & 1&  & \\ \hline
  
  \end{array}$$
  $$\begin{array}{l|cccccccccc} 
 & 1 & x & y & x^2 & xy & y^2 & x^3 & x^2y & xy^2 & y^3\\ \hline
  
   \hbox{introduce } y &  & &  & & & & & & & \\
  (y,0,0) &  &  & 1 &  & & & & & & \\
  (y,y,0) &  &  &  & & 1 & & & 1&  & \\
  (0,y,y) &  &  &  & &  & & &1 & 1 &  \\ \hline
    \hbox{after reduction} &  & &  & & & & & & & \\
    (y,0,0) &  &  & 1 &  & & & & & & \\
  (y+x,y+1,x+1) &  &  &  & &  &  & & &  & \\
  (1,y+1,x+y) &  &  &  & &  & & & & 1 &  \\ \hline
  \end{array}$$
  Hence $(y+x,y+1,x+1)$  is a syzygy and therefore $y+x\in I$
  and it is the first element in deglex. order; we can now omit all the
  multiples of  $y\cdot(1,1,0)$ from consideration. Continuing the
  computation we find
  $$\begin{array}{l|cccccccccccc}
 & 1 & x & y & x^2 & xy & y^2 & x^3 & x^2y & xy^2 & y^3 & x^4 & x^3y\\ \hline
  \hbox{introduce } x^2 &  & &  & & & & & & & &  &\\
  (x^2+x,0,x) &  & &  & & & &  &1 & & & &\\
  (x,x^2+x,0)  &  &  &  & & & & & & & & 1 & \\
  (x,x,x^2) &  &  &  & &  & & & &  & &  & 1 \\ \hline
   \hbox{after reduction }  &  & &  & & & & & & & &  &\\
  (x^2+x+1,1,0) &  & &  & & & &  & & & & &\\
  (x,x^2+x,0)  &  &  &  & & & & & & & & 1 & \\
  (x,x,x^2) &  &  &  & &  & & & &  & &  & 1 \\ \hline
  \end{array}$$
Thus $(x^2+x+1,1,0)$ is a syzygy and $x^2+x+1$ is the second basis
element in $I$ relative to deglex. We can omit all multiples of $
x^2 (1,0,0)$. It follows that $\{ y+x,x^2+x+1 \}$ is the required
Gr\"obner basis.
\end{example}

Note that the above procedure is completely general and can be
used for any base field. Although the general construction uses
only straightforward linear algebra it has a major drawback in
that to determine that a polynomial $f$ belongs to the ideal (in
which case $f$ will be an element of the Gr\"obner basis), one
must compute the minimal representation $f=\sum h_i f_i$ where the
$f_i$ are the initial generators. It is known that the degrees of
the $f_i$ can be doubly exponential in $n$, the number of
variables. This is usually called the Nullstellensatz problem
\cite{Null}.

\begin{remark}\label{r:razones}
However, the particular properties of our setting allow us to use
this algorithm for computing the Gr\"obner basis associated to a
binary code:
\begin{enumerate}
\item \label{ite1}Since the words in our initial generating set
are of the form $w_i-1$, after the first reduction we always have
only elements in $[X]$ as the representative elements for canonical
forms (i.e. coordinate vectors in the vector space $K[X]$). \item
\label{ite3} Because of (\ref{ite1}.) above, in our case, a vector of
$K[X]^{r+1}$ introduced a row either reduces to zero or else it
represents a new irreducible element. Therefore, an element does
not reduce to one of lesser degree apart from to degree zero (in
which case we have obtained a new syzygy and a new element of the
reduced basis). \item \label{ite4} We use a total degree
compatible ordering $<$ on $[X]$ and the new ordering in the
module is a TOP ordering, which looks first for the maximal terms
in any position, and after that takes into account that $\mathbf{e}_i<\mathbf{e}_j$
if $i<j$. \item \label{ite5} From (\ref{ite3}.) and (\ref{ite4}.) above we
find that the degrees of all components in the vectors are the
same, which implies that the degrees of the cofactors (the $h_i$)
are at most the degree of the new element $g$ of the basis. For
this element $g$, the leading term $\T(g)$ is in standard form
(otherwise it would be a multiple of some $x_i^2$  which
contradicts $g\in G_\T\setminus\{x_i^2-1\mid i=1,\ldots, n\}$).
The maximal length of a standard form is $n$.
\end{enumerate}
\end{remark}

\begin{remark}\label{rem:tpat}
Note that since the terms are added in the ordering used for
computing the Gr\"obner basis associated to the code then the
first syzygy we find so that it corresponds to a binomial $g$
whose maximal term is in standard form, satisfies $t=\Td(g)-1$
(see Remark~\ref{rem:t}).
\end{remark}

\begin{example} Consider as a ``toy example" the binary code $\mathcal C$ with generator matrix
$$G=\left( \begin{array}{ccc}
1 & 0 & 1\\
0 & 1 & 1
\end{array}\right).$$
We find that \begin{equation*}\begin{split} I(\mathcal{C})= &
\langle f_1=x_1x_3-1,
f_2=x_2x_3-1, \\
& f_3=x_1^2-1,  f_4=x_2^2-1, f_5=x_3^2-1 \rangle.\end{split}
\end{equation*}
In the associated syzygy computation the rows corresponding to the
binomials $x_i^2-1$ are considered as implicit in the
computations: see, for example, in the Table below when the syzygy
corresponding to $x_3-x_1$ is obtained.

$$\begin{array}{l|cccccccc|c}
 & -1 & x_1 & x_2 & x_3  & x_1x_2 & x_1x_3 & x_2x_3 & x_1x_2x_3 & \hbox{multiples of } x_i^2\\ \hline
 (1,0,0,0,0,0) & 1 &  &  &  & & & & & \\
 (1,1,0,0,0,0) &  &  &  &  & & 1& & & \\
 (1,0,1,0,0,0) &  &  &  &  & &  & 1 & & \\ \hline
 \hbox{introduce } x_1 &  &  &  &  & &  &  & & \\
 (x_1,0,0,0,0,0) &  & 1 &  &  & & & & & \\
 (x_1,x_1,0,0,0,0) &  &  &  &  & & & & & x_1^2x_3 \\
 (x_1,0,x_1,0,0,0) &  &  &  &  & &  &  & 1& \\
  \hline
  
 \hbox{introduce } x_2 &  &  &  &  & &  &  & & \\
 (x_2,0,0,0,0) &  &  & 1  &  & & & & & \\
 (x_2,x_2,0,0,0,0) &  &  &  &  & & & & 1 &  \\
 (x_2,0,x_2,0,0,0) &  &  &  &  & &  &  & & x_2^2x_3 \\
  \hline
 \hbox{reduction }  &  &  &  &  & &  &  & & \\
 (x_2,0,0,0,0,0) &  &  & 1  &  & & & & & \\
 (x_2-x_1,x_2,x_1,0,0,0) &  &  &  &  & & & &  &  \\
 (x_2,0,x_2,0,0,0) &  &  &  &  & &  &  & & x_2^2x_3 \\
  \end{array}$$
Thus $x_2-x_1=x_1f_2-x_2f_1$, $x_2-x_1$ belongs to the Gr\"obner
basis and we can now omit all the multiples of
$x_2(1,1,0,0,0,0,0)$ from our computation. Continuing we find
$$\begin{array}{l|cccccccc|c}
 & 1 & x_1 & x_2 & x_3  & x_1x_2 & x_1x_3 & x_2x_3 & x_1x_2x_3 & \hbox{multiples of } x_i^2\\ \hline
 \hbox{introduce } x_3 &  &  &  &  & &  &  & & \\
 (x_3,0,0,0,0,0) &  &  &   &1  & & & & & \\
 (x_3,x_3,0,0,0,0) &  &  &  &  & & & &  & x_1x_3^2 \\
 (x_3,0,x_3,0,0,0) &  &  &  &  & &  &  & & x_2x_3^2 \\  \hline
 \hbox{reduction }  &  &  &  &  & &  &  & & \\
 (x_3,0,0,0,0,0) &  &  &   &1  & & & & & \\
 (x_3-x_1,x_3,0,0,0,x_1) &  &  &  &  & & & &  &  \\
 (x_3-x_2,0,x_3,0,0,x_2) &  &  &  &  & &  &  & & \\  \hline
  \end{array}$$
We have the syzygy  $x_3-x_1=x_1f_5-x_3f_1$ and $x_3-x_1$ belongs
to the Gr\"obner basis (note that we can make reductions of terms
$T\cdot x_i^2$, $T$ a term, as soon as we have introduced $T$
since $x_i^2-1$ is a generator). The result of the computation is
the Gr\"obner basis $\{x_2-x_1, x_3-x_1, x_1^2-1\}$

\begin{remark}
In recording the computations we need only to keep the first
components on the left and pointers to those places with a 1 in
the rest of the table. This gives the following adapted FGLM basis
conversion algorithm.
\end{remark}
\end{example}

\subsection{Adapted FGLM algorithm}\label{sb:fglm-adap}
The algorithm  computes a Gr\"obner basis for the syzygy module,
but we are interested only in the first component which is the
Gr\"obner basis $G_\T$ for $I({\mathcal C})$. For theoretical
reasons we will denote by $G(M)$ the set which is constructed by
the algorithm, which on termination is a Gr\"obner basis for the
module, but we will just compute the first component $G$ of this
set.

In the algorithm we use two main structures. One is {\em List},
which has the form $(v_1,v_2)$, where $v_1$ represents the first
component of the corresponding vector in the module, and $v_2$ is
the representative element in $K[X]$ (in our case an element of
$[X]$ -- see (\ref{ite1}) of Remark~\ref{r:razones}). If
$w=(v_1,v_2)\in \mathit{List}$ then we write $w[1]=v_1$ and
$w[2]=v_2$. The second structure is the list $N$ that stores the
first components of the elements of {\em List} that are canonical
forms. The third structure is the list $V$ whose  $r$-th element $v_r$  is the representative
element in $K[X]$ of the $r$-th element of $N$ (the second component of the pairs in {\em List}).

\noindent Subroutines of the algorithm:
\begin{itemize}
\item ${\bf InsertNexts}(w, \mathit{List})$ inserts the products
$wx$ (for $x \in X$) in {\em List} and sorts it by increasing
order with respect to $<$, with account being taken first of the
first component, and, in case these are equal, then by comparison
of the second components. The reader should note that ${\bf
InsertNexts}$ could count the number of times that an element $w$
is inserted in {\em List}, so $w[1] \in N_<(I) \cup T_<\{G\}$ if
and only if this coincides with the number of variables in the
support of $w[1]$ (if not, this would means that $w[1]\in
T_<(I)\setminus T_<\{G\}$, see \cite{FGLM}). This criterion can be
used to determine the boolean value of the test condition in
Step~4 of the Algorithm~\ref{a:fglm-syzygy} \item ${\bf
NextTerm}(\mathit{List})$ removes the first element from {\em
List} and returns it. \item \textbf{Member}($v,[ v_1,\ldots
,v_r]$) returns $j$ if $v=v_j$ or false otherwise.
\end{itemize}

\begin{alg}\label{a:fglm-syzygy} $ $
\begin{description}
\item[\bf Input] $F=\left\lbrace w_1-1,w_2-1,\ldots
,w_r-1\right\rbrace $ the set of
binomials associated with a generating set of a binary code\\
$<_\T$ a total degree compatible ordering 
\item[\bf Output] The
reduced Gr\"obner basis $G_\T$ of the ideal \\ $\left\langle
F\cup\left\lbrace x_i^2-1\mid i=1,\ldots, n\right\rbrace
\right\rangle$ w.r.t. $<_\T$
\end{description}
\begin{enumerate}
\item $\mathit{List}:=[(1,1),(1,w_i)_{i=1,\ldots,r},
(1,x_i^2)_{i=1,\ldots,n}]$ (the elements should be ordered
following $<_\T$ in the second component of the pairs),\\
$G_\T:=\{~\}$,$N:=[~]$
%$\mathrm{Mult}=\{\}$.
\item While $\mathit{List}\neq \emptyset$ do \item  $\qquad
w:=\mathrm{NextTerm}(\mathit{List})$; \item $\quad\;\;$  If $w
\notin \T(G(M))$; \item  $\qquad v^\prime:=w[2]$; \item  $\qquad
j:= \mathrm{Member}(v^\prime ,[v_1,\dots ,v_r])]$; \item
$\qquad$If $j \neq$ false $\;$then $G:=G \cup \{w[1]-w_j\}$;
%\item[8.] $\qquad$ $\qquad\qquad \phi(u,x_k):=w_j$
\item $\qquad$ else $r:=r+1$; \item $\qquad\qquad\,
v_r:=v^\prime$; \item  $\qquad\qquad\, w_r:=w[1],\; N:=N\cup
\{w_r\}$; \item  $\qquad\qquad
\mathit{List}:=\mathrm{InsertNexts}(w_r,\mathit{List})$; \item
Return[$G$]
%\end{center}
\end{enumerate}
\end{alg}

Note that this algorithm for computing the Gr\"obner basis  $G_\T$ associated to the code $\mathcal{C}$ is especially well suited in our setting since all the elements in the basis (respectively codewords, cycles) appear in an increasing term ordering (respectively increasing ordering on the weight or the length) during the computation. Moreover,  the computation can be stopped when a desired weigth of the codewords (respectively length of the cycles) is obtained which is usefull  for finding  many  combinatorial properties of the code (respectively the graph) such that the minimal distance (see Remarks \ref{rem:t}, \ref{rem:tpat}) or finding the minimal codewords (see Proposition \ref{p:minpes}).

\end{document}